\documentclass[10pt]{amsart}%
\usepackage{amsmath}
\usepackage[dvips]{graphicx}%
\usepackage{amsfonts}%
\usepackage{amssymb}
\newtheorem{theorem}{Theorem}[section]
\theoremstyle{plain}

\newtheorem{corollary}[theorem]{Corollary}

\newtheorem{lemma}[theorem]{Lemma}

\newtheorem{proposition}[theorem]{Proposition}

\theoremstyle{definition}

\newtheorem{remark}[theorem]{Remark}
\theoremstyle{remark}
\numberwithin{equation}{section}
\newcommand{\inclupdef}{\def\inclup{\lower2pt\vbox{\lineskiplimit\maxdimen
\lineskip-3.25pt\ialign{$\mathsurround=0pt####$\crcr\kern.625pt\vrule height5pt
width.2pt\hfil\vrule height5pt
width.2pt\kern.625pt\crcr\smallsmile\crcr}}\kern1.5pt\vrule height6.625pt
width0.2pt depth0.625pt}}
\inclupdef

\renewcommand{\geq}{\geqslant}

\begin{document}
\title[ Irreducible subfactors ]{Irreducible subfactors
derived from Popa's construction for non-tracial states}
\author{Florin R\u{a}dulescu}
\address{Department of Mathematics, The University of Iowa, Iowa 
City, IA 52242-
1419 U.S.A.}
\email{radulesc@math.uiowa.edu}
\urladdr{http://www.math.uiowa.edu/\symbol{126}radulesc/}

\begin{abstract}
For an inclusion of the form
  $\Bbb C\subseteq M_n(\Bbb C)$, where $M_n(\Bbb C)$ is endowed with a 
state with diagonal
weights $\lambda=(\lambda_1, ..., \lambda_n)$,  we use  Popa's 
construction, for
non-tracial states, to obtain  an
  irreducible  inclusion of $II_1$ factors, $N^\lambda(Q)\subseteq 
M^\lambda(Q) $ of index $\sum
\frac{1}{\lambda_i}$.
$M^\lambda(Q)$ is identified with  a subfactor inside
   the centralizer algebra of the canonical free product state on 
$Q\star M_N(\Bbb C)$. Its
structure is described by  ``infinite'' semicircular elements as in \cite{Ra3}.

The irreducible subfactor inclusions obtained by this method are 
similar to the first
irreducible subfactor inclusions, of index in $[4,\infty)$ constructed in
\cite {Po1}, starting with the
  Jones' subfactors inclusion $R^s\subseteq R$, $s>4$. In the present 
paper, since
the inclusion we start with has a simpler structure,  it is  easier 
to control the algebra
structure of the subfactor inclusions.

  If the weights correspond
  to  a unitary, finite dimensional representation
of a Woronowicz's compact quantum group $G$, then the  factor $M^\lambda(Q)$ is
contained in the fixed point algebra of an action of the quantum 
group on $Q\star M_N(\Bbb C)$ , with equality
if $G$ is $SU_q(N)$, (or   $SO_q(3)$ when $N=2$).
  By Takesaki duality,   the factor $M^\lambda(L(F_N))$  is Morita
equivalent to $\mathcal L(F_\infty)$.

This
method gives   also another approach to find,
as also recently proved in \cite {ShUe},
  irreducible subfactors of $\mathcal L(F_\infty)$ for index values 
bigger than 4.
\end{abstract}
\maketitle

\setcounter{section}{-1}
\section{\label{Int}Introduction and definitions} In this paper
  we consider the structure of  subfactors obtained from Popa's 
construction, for non-tracial
states,  for the  inclusion $M_N(\Bbb C)\subseteq M_N(\Bbb C)\otimes 
M_N(\Bbb C)$.
We  fix a diagonal matrix with non-zero weights
$\lambda_1,...\lambda_N$. The states on the two algebras
  are respectively $tr (D\cdot)$ and $tr (D\otimes D^{-1}\cdot)$. This 
is then the
Jones' iterated basic construction (\cite {Jo}, \cite {Jo1})
  for $\Bbb C\subseteq M_N(\Bbb C)$, where $M_N(\Bbb C)$ is endowed 
with the state $tr
(D\cdot)$. The algebra
$M_N(\Bbb C)\otimes 1$ is invariated by the modular
  group of the state $tr (D\otimes D^{-1}\cdot)$ on $M_N(\Bbb 
C)\otimes M_N(\Bbb C)$. Hence by
\cite{NaTa}, there exists a conditional expectation
  and a corresponding  Jones projection, which is also invariated by
  the modular group. The value
of the state on the Jones projection is $\sum \frac{1}{\lambda_i}$.
  Let $Q$ be a diffuse finite von Neumann algebra with faithful trace $\tau$.

We
apply Popa's construction for the inclusion
$M_N(\Bbb C)\subseteq M_N(\Bbb C)\otimes M_N(\Bbb C)$.
  This yields an irreducible inclusion $N^\lambda(Q)\subseteq 
M^\lambda(Q)$ of type
$II_1$ factors, inside the type $III$ factor
  $(Q\otimes M_N(\Bbb C)) *_{M_N(\Bbb C)} (M_N(\Bbb C)\otimes M_N(\Bbb C))$,
  considered with the
amalgamated free product state.

By applying the reduction method described in \cite{Ra}, \cite{Ra4}, 
we can reduce this
  procedure to
the case when the algebra over which we amalgamate is abelian. Thus, 
we end up with a concrete
description of the factor and the subfactor inside $Q\star M_N(\Bbb 
C)$, where  the later
algebra is endowed with the free product state $\tau\star tr 
(D\cdot)$.  This factor for
$N=2$, $Q=L^\infty([0,1])$ can be very explicitly described in terms 
of the ``infinite''
Voiculescu's type (\cite {Vo}) free semicircular element used to 
identify the core of $
L^\infty([0,1])\star M_2(\Bbb C)$ in \cite{Ra3}.

Moreover, by analogy with the case of a trace on $M_N(\Bbb C)$, it 
turns out that  $M^\lambda(Q)$ is  contained  in the fixed
  point algebra for a
free  product coaction (of the type considered in \cite{Ue1}) on 
$Q\star M_N(\Bbb C)$.

Let $\alpha$ be any finite dimensional unitary
  selfadjoint coaction of a Woronowicz quantum group $G$ on the finite 
dimensional algebra
  $M_N(\Bbb
C)$. We prove that the algebra $M^\lambda(Q)$ described above is 
contained in the fixed point
algebra of coaction of the quantum group. If $G$ is the quantum group 
$SU_q(N)$ then the
algebra $M^\lambda(Q)$ is exactly the fixed point algebra. This is 
analogous to the result
in \cite {Ko}, \cite {Na}, where the fixed point algebra of the 
infinite product action on
$M_N(\Bbb C)$
of $SU_q(N)$ coincides with the hyperfinite factor.

Assume (by analogy with \cite {Wa}) that
$\alpha^{\otimes^n}$ contains any other finite dimensional unitary 
representation of $G$.
  Then we prove that the fixed point algebra
$(Q\star M_N(\Bbb C))^G$ is Morita equivalent to the cross
product by $\alpha$. Such a cross product   is naturally  described 
(\cite{Ue1}) by a  free product with amalgamation.
Since the fixed point
  algebra is a $II_1$ factor it follows that the fixed point algebra 
is Morita equivalent with an amalgamateted free product of the form
  $(Q\otimes D)\ast_D C$, where $C,D$ are (infinite) direct sums of 
matrix algebras. By the techniques in
\cite{Ra} it follows, for
$Q=\mathcal L(F_N)$, that the fixed point algebra is
the von Neumann algebra of a free group with infinitely many generators.
  The $II_1$ factors  is invariated
by the action of an automorphism scaling the trace on a larger 
$II_\infty$ factor.

We also reobtain in this way, by a different approach, the
result recently proved in the remarkable paper by D. Shlyahtenko and 
Y. Ueda, \cite {ShUe}, that $L(F_\infty)$ has irreducible subfactors 
for any
index value bigger than 4.  One theme in their paper is   that one can
determine the isomorphism class of a fixed point algebra by looking 
at the crossed product
algebra (via Takesaki's duality).  We owe to their paper  the use of this
philosophy for  coactions.

We do not know if the subfactors in this paper coincide the ones constructed in
\cite {ShUe}. In both cases the higher relative commutants invariant is
$A_\infty$ (though the factors in \cite {ShUe} seem to generate other 
higher relative
commutants invariants) when used for other compact quantum groups). 
Both algebra of the factors
are obtained by taking the  the fixed point algebra of a coaction of 
a quantum group. In one
case the fixed point algebra is a
$II_1$ factor
  and in the other case the fixed point algebra is of type $III$.

  The construction in this paper  proves  that the non-tracial version of Popa's
construction of subfactors,  naturally yields irreducible subfactors, 
In addition the
corresponding algebras are fixed algebras of coactions of quantum groups.

These algebras, although they are very similar to the subfactors in 
the breakthrough
construction in
\cite{Po1}, are not yet proven to be isomorphic to the algebras in 
\cite {Po1},\cite{Po3}.
Though, our result is strong evidence to the  conjecture that the 
subfactors in \cite{Po1},
\cite{Po2} are free group factors when $Q$ is $\mathcal L(F_\infty)$ 
(the only case when
the subfactors in \cite{Po1}, \cite{Po2} are  known to be free group 
factors is for index
values less than 4, and the higher relative commutants are the 
Temperly-Lieb algebra, see
\cite{Ra}).

  The method in this papers  also gives an explicit model of the 
subfactors in terms of
(infinite) semicircular random matrices (\cite{Vo},\cite {Ra2},\cite {Ra3}).

This work was partially  by NSF
grant $DMS$.  The author is  grateful
to T. Banica, D. Bisch, F. Boca, F. Goodman, G. Nagy, M. Pimsner,
   D. Shlyahenko,V. Toledano,  Y. Ueda for  disscutions
on this subject. This work was initiated during the author's stay at 
the University of Geneva,
and was partially supported by the Swiss Science National 
Foundation. The author is very
thankful to P. de la Harpe for the  very warm welcome at the 
University of Geneva.

\section{\label{De1} Subfactors derived from Popa's construction for 
non-tracial states }

In this section we  describe the structure
  of the algebra of the subfactors that are
  obtained, by using Popa's construction,
from a finite dimensional inclusion $B\subseteq A$
  with Markov state. We  apply this construction to the case
of the inclusion $\Bbb C\subseteq M_n(\Bbb C)$,
  where $M_n(\Bbb C)$ is endowed with a state. This will give
  irreducible subfactor of indices from
$4$ to $\infty$.

  We start with the following lemma which proves
  that the subfactor associated to $B \subseteq A$ only depends on the 
inclusion matrix. To do
this we use a reduction style procedure that was used in \cite{Ra}, 
to the case when $B$ is
abelian. The following lemma is proved in \cite{Ra4}

\begin{lemma}
\label{LemDe.1} Let $C\subseteq B$ be a
  finite dimensional inclusion of matrices, with trivial centers 
intersection and assume that
$B$ is endowed with a $\lambda$ Markov
  state (that is there exists a normal
  conditional expectation from $B$ onto $C$ which
preserves the state, there exists a state  on the corresponding basic 
construction extending
the given one  and the corresponding Jones' projection has expectation
$\lambda$ times the scalars).

   Let $Q$ be a type $II_1$ factor, let $A=<B,e_1>$ be the first step 
in the Jones' basic
construction of
$C\subseteq B$, where  $e_1$ is the Jones projection. Assume that $A$ 
is endowed with the canonical Markov state.
Let $m_i\in B$ be a Pimsner-Popa
orthonormal basis for the inclusion $C\subseteq B$.
  Let $\mathcal A$ be the algebra  $(Q\otimes B) \star_B A$ and let 
$\Phi$ be the Popa's map
associated to the inclusion defined by
$$\Phi (x)=\sum m_ie_1 x e_1m_i^*, x\in \mathcal A.$$
Then as proved in \cite{Po1},\cite{Po2},\cite{Po3}, $\Phi$ maps $B'$ into $A'$.

Let $f_i$ be a family of representatives of minimal
  projections in $B$ and let $F_0$ be the sum of this projections.
Let $\mathcal A_0, A_0, B_0$ be the reductions
  of respectively the algebras $\mathcal A, A, B$ by the projection 
$F_0$. In particular
$B_0$ is abelian and, as proved
  in (\cite{Ra}), $\mathcal A_0$ is identified with  $(Q\otimes 
B_0)\star_{B_0} A_0$.
Since $F_0$ has central support 1, $\Phi$
  induces a map on $\mathcal A_0$ which maps $B'_0 $ into $A'_0$.
Let $N(Q)\subseteq M(Q)\subseteq B'$ be the
  minimal algebra (introduced in \cite{Po1}) in $\mathcal A$ that is 
closed under $\Phi$ and
contains
$Q$.  Correspondingly let $N_0(Q)\subseteq M_0(Q)\subseteq B_0'$
  be the minimal algebra in $\mathcal A_0$
that is closed under $\Phi_0$ and contains $Q$.
Then
$$ M_0(Q)=(M(Q))_{F_0}; N_0(Q)=(N(Q))_{F_0}. $$
Moreover if $m'_i$ is any orthonormal basis for $B_0\subseteq A_0$, then
$$\Phi_0(x)=\sum m'_ix(m'_i)^*, x \in \mathcal A_0.$$
\end{lemma}
\begin{proof} This lemma will be proved in full generality in \cite 
{Ra4}. We  give here
in this section an ad-hoc proof for the case $C=\Bbb C\subseteq B= 
M_N(\Bbb C)$, where
$M_N(\Bbb C)$ has a state $tr(D\cdot)$ (as in the next lemma). Indeed 
the next step in the
basic construction (as proved in the next lemma) is $A=M_N(\Bbb 
C)\otimes M_N(\Bbb C)$.
Here the algebra $\mathcal A$ is
$(Q\otimes M_N(\Bbb C))\ast_{M_N(\Bbb C)} (M_N(\Bbb C)\otimes M_N(\Bbb C))$ and
consequently $\mathcal A_0$ is $Q\ast M_N(\Bbb C)$.
By construction
$\Phi_0$ is of the form $\Phi_0(x)=\sum_\alpha n_\alpha x 
(n_\alpha)^*$, for some
fixed $n_\alpha$ in $M_N(\Bbb C)$. Let $e_{ij}$ be a matrix unit 
diagonalizing $D$.

By cutting by a minimal projection, it follows that $\Phi_0(x)$ takes 
values into
$(M_N(\Bbb C))' \cap \mathcal A_0$.  Thus necessary, there exists 
real $\theta_1,...,\theta_N$
numbers such that
$\Phi_0$ is of the form
$$\Phi_0(x)=\sum \theta_j e_{ij}xe_{ji}, x\in \mathcal A_0.$$
But $\Phi_0$ has also the property that the conditional expectation from
$\mathcal A_0$ onto $M(Q)$ maps $M_N(\Bbb C)$ into the scalars. This is only
possible if $\theta_j e_{ij}$ is an orthonormal basis for $M_N(\Bbb C)$, i.e.
if $\theta_i$ are the inverses of the eigenvalues of $D$.
\end{proof}

  As  will see bellow, even in the case
when the inclusion $B\subseteq A$
  is $\Bbb C\subseteq M_N(\Bbb C)$, with a trace on $M_N(\Bbb C)$
   has to be handled by a rather
a complicated machinery.
  Indeed in this case the above subfactor  is the fixed algebra under
the action of $SU(N)$ on $Q\ast M_N(\Bbb C)$ (here $SU(N)$ acts 
trivially on Q and by
conjugation on $M_N(\Bbb C))$.
  This is because the element $\sum e_{ij}\otimes e_{ji}$ is
invariant under the product action of the group $SU(N)$  and all the 
others are obtained by
intercalation or concatenation from this element (\cite {Wo3}).

We generalize bellow the above construction to the case when
  the trace on $M_N(\Bbb C)$ is replaced by a state.
First we record the following well known  folklore lemma.
  It deals with the Jones basic construction for states (see e.g. 
\cite{Po4}, \cite{Jo1})
\begin{lemma}
\label{LemDe.3} Let $\phi=tr(D\cdot) $ be
  the state on $M_n(\Bbb C)$ where $D$ has
  diagonal eigenvalues $\lambda_1,...,\lambda_n$.
Then the next step in the Jones basic construction
  with Markov state is $M_n(\Bbb C)\otimes M_n(\Bbb C)$ with state
$$\frac{1}{tr D^{-1}}tr ((D\otimes D^{-1}\cdot).$$
  This state is so that the the Jones projection is invariated by the 
modular group and
projects  onto a scalar multiple of the identity. The value of the 
scalar is $(\sum
\frac{1}{\lambda_i})^{-1}$.

\end{lemma}
\begin {proof}

Indeed if $(e_{ij})$ is a matrix unit, which also diagonalizes $D$, then
the Jones' projection is

$$\sum _{ij} \lambda_i^{-1\slash 2}\lambda_j^{-1\slash 
2}e_{ij}\otimes e_{ij}.$$
\end{proof}

In the next lemma we describe what  Popa's construction (\cite{Po1}) is in
  the non-tracial case in a very particular case (for the more general 
construction
see \cite{Ra4}). This
construction allows
  to obtain irreducible subfactors,
  in factors with non-trivial fundamental group,
of index values in
$[4,\infty)$ starting from the
  inclusion $\Bbb C\subseteq M_N(\Bbb C)$.
  We will describe this factors explicitly.

\begin{lemma}
\label{Popant}
  We use the notations from the previous lemma and its proof.
  Let $e_1$ be the Jones projection for the basic construction
  in the previous lemma.
Let $(M_\alpha) e_1$, in the centralizer algebra of the state $\phi$,
  be a Pimsner-Popa basis for
the inclusion
$M_N(\Bbb C)
\subseteq M_N(\Bbb C)
\otimes M_N(\Bbb C).$ Note that this is allways possible since the centralizer
algebra (which contains $e_{ij}\otimes e_{ij}$), is isomorphic to 
$M_N(\Bbb C)$ and contains
$e_1$.

Let $Q$ be a diffuse abelian von
  Neumann algebra or a type $II_1$ factor with faithful trace
$\tau$
and let
  $\mathcal A=(Q\otimes M_N(\Bbb C))\ast_{M_N(\Bbb C)} (M_N(\Bbb C)
\otimes M_N(\Bbb C))$ be the amalgamated free product factor obtained 
by taking the
GNS-construction corresponding to the given state on $M_N(\Bbb C)
\otimes M_N(\Bbb C)$ and the trace $\tau$. This is obviously a type 
$III$ factor if
$\phi$ is not a trace (\cite{Dy},\cite{Shly1}, \cite{BlDy},\cite{Ue3}).

Let as in \cite {Po1},\cite{Po2}, $\Phi$ be the map defined on
$\mathcal A$ by
$$\Phi(x)=\sum_\alpha M_\alpha e_1 x e_1 M_\alpha^*, x\in \mathcal A$$
Let $M^\lambda(Q)$ be the minimal algebra in $\mathcal A$ containing 
$Q$ that is invariant
under the map $\Phi$. Let $N^\lambda(Q)$ be the image of $M^\lambda(Q)$ through
$\Phi$.

  Then, (as in \cite {Po1})
$N^\lambda(Q)\subseteq M^\lambda(Q)$ is an irreducible inclusion of 
type $II_1$ factors
of index $\sum \frac{1}{\lambda_i}$.
\end{lemma}
\begin{proof} This is basically proved in \cite {Po1},\cite {Po2} (see also
\cite{Ra4}). The only
additional care comes from the fact that we are dealing with a state 
instead of a trace.
But since $M_\alpha e_1$ are chosen in the centralizer of the state 
it follows that
$M^\lambda(Q)$ is a type $II_1$ factor. The main part of the argument 
(the Markovianity of
the trace) follows from the fact that the conditional expectation of 
the algebra
$M_N(\Bbb C) \otimes M_N(\Bbb C)$ onto $M^\lambda(Q)$ is the algebra 
of scalars,
which follows from the properties of the Pimsner Popa basis.

\end{proof}
It is not clear if the above factor is the same as the one 
constructed in \cite {Po1}
starting from the inclusion derived from the Temperly-Lieb algebra 
for index values
bigger than 4. In that case the algebras stay in a much larger amalgamated free
product, since the Pimsner-Popa basis sits in hyperfinite factor (and 
can not be chosen
to belong to a finite dimensional algebra).
\begin{remark}
\label{addsteps}
  To obtain any additional step in the Jones' basic construction of 
the above inclusion
one has to proceed as follows. One starts with the iterated Jones's 
basic construction
of $\Bbb C\subseteq M_N(\Bbb C)$, where $M_N(\Bbb C)$ is endowed as 
above with the
state $tr(D\cdot)$. Then the  $k$-th, $(k+1)$-th steps of this basic 
construction are
$B=M_N(\Bbb C)^{\otimes^k}$, $A=M_N(\Bbb C)^{\otimes^{k+1}}$. The 
states on the iterated steps
in the basic construction are states with weights involving 
consecutive products of
$D$ and $D^{-1}$. Let $e_{k-1},...,e_1$be the corresponding Jones'
projections, indexed so that $e_1\in A$ is exactly the Jones projection
for $M_N(\Bbb C)^{\otimes^{k-1}}\subseteq B.$

    Then we perform the same construction as
above for $B\subseteq A$ (instead of $M_N(\Bbb C)\subseteq M_N(\Bbb 
C)\otimes M_N(\Bbb C)$).
We get a subfactor inclusion $N^\lambda(Q)\subseteq M^\lambda(Q)$. 
The iterated Jones' basic
construction steps (up to $k-1$) are then obtained by adding 
consecutively the projections
$e_1,e_2,..., e_{k-1}$. As we will see bellow by reducing by a 
minimal projection in $B$
it follows that these subfactors are isomorphic to the original ones. 
In particular
it follows that the factor and subfactor are allways isomorphic to 
the algebra corresponding to
$D$ and
$D^{-1}$
\end{remark}

\begin{lemma}
\label{LemDe.4} Let $\phi $ be a
  state on $M_N(\Bbb C)$ with $\phi=tr(D.)$,
  where $D$ is diagonal with eigenvalues list
$\lambda_1,..., \lambda_N$. Let $(e_{ij})$ be
  a matrix unit for $M_N(\Bbb C)$, such that
  $e_{ii}$ are the projections onto the eigenvectors
of
$D$. Let $Q$ be a type $II_1$ factor and
  consider the type $III$ factor $Q\ast M_N(\Bbb C)$ (see e.g 
\cite{Bar}, \cite{Dy1},
\cite{Ra3},
\cite{Shly1})
  where the free product is with respect to the state
  $\phi$ on $M_N(\Bbb C)$ and the trace on $Q$.
Let $$\Phi(\alpha)=
\sum_{ij} \lambda_j^{-1} e_{ij}\alpha e_{ji},\ \ \ \ \ \alpha \in 
Q\ast M_N(\Bbb C),$$
and let $\mathcal M^\lambda(Q)$ be the minimal
  subalgebra of $Q\ast M_N(\Bbb C)$ that contains $Q$ and is closed 
under $\Phi$.
Let $\mathcal N^\lambda(Q)$ be the image of $M$ through $\Phi$. Then
$$[\mathcal M^\lambda(Q):\mathcal N^\lambda(Q)]=
\sum_i\frac{1}{\lambda_i} ;\ \ \
   \mathcal M^\lambda(Q)\cap\mathcal N^\lambda(Q)'=\Bbb C1.$$
Moreover $\mathcal M^\lambda(Q)$ has non-trivial fundamental group and
  the type of the subfactor inclusion is $A_\infty$. The algebra
of the subfactor is isomorphic with the factor
$M^{(\lambda_i^{-1})}(Q)$.
\end{lemma}
\begin{proof} We apply Popa's construction
  for the inclusion $M_N(\Bbb C)\subseteq M_N(\Bbb C) \otimes M_N(\Bbb 
C) $ described in
Lemma
\ref{LemDe.3}. Since this is the basic
  construction for $\Bbb C\subseteq M_N(\Bbb C)$ we can apply
Lemma
\ref{LemDe.1} and one obtains an
  equivalent description of the subfactor $\mathcal M^\lambda(Q)$ sitting
  in $Q\ast M_N(\Bbb C)$. Since $\lambda_j^{-1} e_{ij}$ is a Pimsner-Popa
basis for the inclusion
$\Bbb C\subseteq M_N(\Bbb C)$, the result follows. That the subfactor 
inclusion is of type $A_\infty$ is proved in a more general context in
\cite{Ra4}.
\end{proof}

By analogy with the case of traces,  we have the following result 
which shows that the algebra obtained by Popa's construction from
the inclusion of finite dimensional algebras (with states), $\Bbb 
C\subseteq M_N(\Bbb C)$ is a minimal fixed point algebra.

\begin{theorem}
\label{TheDe.5} Let $A$ be the function algebra of Woronowicz quantum 
group $G$  and let $V_\pi$ be finite
dimensional unitary representation of dimension $N\geq2$.
Assume that the operator $Q_\pi$ associated to the representation as 
in \cite{Wo}, \cite{Ba} has eigenvalues
$\lambda_1^{-1},...,\lambda_N^{-1}$.  Let $Q\ast B(V_\pi)$ and 
consider the Ueda's
(\cite{Ue1},\cite{Ue2}) type  of free product coaction on $Q\ast 
B(V_\pi)$, that acts
trivially on $Q$ and by conjugation on $B(V_\pi)$. Let $\mathcal M$ 
be the Popa's type  factor
constructed in lemma \ref{LemDe.4}.

Then
$$\mathcal M\subseteq(Q\ast B(V_\pi))^G.$$
If $A$ is the function algebra of $SU_q(2)$, and $\pi$ is the 
fundamental representation then
we have equality.

\end{theorem}

\begin{proof}Let $(e_{st})$ be a matrix unit for $B(V_\pi)$ that 
diagonalizes $Q_\pi$. Assume that the unitary implementing the
representation is represented as
$$U=\sum_{s,t} e_{st}\otimes u_{st}.$$
Note that we use a another indexing for the entries of the unitary 
than the one used in \cite{Wo}, section 4. Since as proved
in \cite{Ba},
$$Q^{1\slash 2}_\pi {\overline U} (Q_\pi)^{-1\slash 2}$$
is the unitary corresponding to the conjugate representation it follows
that the matrix
$$\alpha_{ij}=\lambda_i^{-1\slash 2}\lambda_j^{1\slash 2}u_{ji},$$
is a unitary too.
Thus we have the following unitarity conditions :
\begin{equation}
\sum_{j}u_{sj}u^*_{tj} =\delta_{st}; \sum_j u^*_{js}u_{jt}= 
\delta_{st}, s,t=1,...N,
  \label{eqDe.4}
\end{equation}

\begin{equation}
\sum _j\lambda_j^{-1} u^*_{sj} u_{tj}=\lambda_s^{-1}\delta_{st};\ \ \ \ \
\sum_{j}\lambda_j u_{js}u^*_{jt}=\lambda_s\delta_{st}; s,t=1,...N.
  \label{eqDe.5}
\end{equation}

The corepresentation for the algebra $B(V_\pi)$ obtained by 
conjugation by the unitary $U$ is given by
\begin{equation}
U (e_{ij}\otimes 1)U^*=\sum_{rs} e_{rs}\otimes u_{ri} u^*_{sj}
  \label{eqDe.6}
\end{equation}
The elements in the algebra $\mathcal M$ have an open and closing 
parenthesis structure
(\cite {Po1},\cite{Bo1}). This proves that these elements are fixed
points under the coaction, by recursively using  (first equality in 
each of) the  relations
\ref{eqDe.4}, \ref{eqDe.5}.
For example
the coaction on an element of the form
$$x=\sum_{\alpha,\beta} \lambda_\beta^{-1}q_1e_{\alpha,\beta}q_2 
e_{\beta,\alpha} q_3$$
gives
\begin{equation}
x=\sum_{r_1,r_2,s_1,s_2} q_1e_{r_1,s_1}q_2e_{r_2, s_2}q_3(\sum_{\alpha,\beta}
\lambda_\alpha^{-1}u_{r_1\alpha}u^*_{s_1\beta}u_{r_2\beta}u^*_{s_2\alpha}).
\label{eqDe.8}
\end{equation}
By applying first equation \ref{eqDe.4}  and then \ref{eqDe.5} it follows that
\begin {equation}
\label{marea}
\sum_{\alpha,\beta}
\lambda_\alpha^{-1}u_{r_1\alpha}u^*_{s_1\beta}u_{r_2\beta}u^*_{s_2\alpha}=\lambda^{-1}_{s_1}
\delta_{r_1s_2}\delta_{s_1r_2}.
\end{equation}
This proves that the coaction on $x$ in (\ref{eqDe.8}) gives 
$x\otimes 1$. The only other fixed point are obtained by recurrence. 
Note that
the algebra $Q$ is in the fixed point algebra of the coaction.

More precisely if $x$ in the fixed point algebra and $h$ is Haar 
measure, we need to show that
$$(Id\otimes h)(\sum_{\alpha\beta} 
\lambda_\beta^{-1}(U(e_{\alpha,\beta}\otimes 1)U^*)
(x\otimes1)(U (e_{\beta,\alpha} \otimes 1) U^*))=
\sum_{\alpha\beta} \lambda_\beta^{-1}e_{\alpha,\beta}xe_{\beta,\alpha}  .$$
This follows by using the modular properties of the Haar measure and 
using the relations
\ref{eqDe.8}. Once we have shown this and since $Q$ is obviously in 
the fixed point algebra
the result follows from the definition of the algebra $M^\lambda(Q)$ 
as the minimal algebra
containing $Q$ closed to the operation in the statement of lemma \ref{LemDe.4}.

Let $Q_0$ be the subspace of $Q$ consisting of elements of zero trace.
The fixed point algebra and Popa's algebra admit a filtration given 
by the subspaces\break
$B(V_\pi)Q_0 B(V_\pi)... Q_0 B(V_\pi) $
  corresponding to the number of
occurrences of copies of $B(V_\pi)$. Moreover the averaging argument 
(with respect to
Haar measure) used in
\cite {Nag} or
\cite{Ko} to establish that the fixed point algebra in the Powers 
factor is the algebra
generated by the Jones projections, allows
to reduce the determination of the the fixed point algebra, to the 
determination of the
intersection of the fixed point algebra with the subspaces in the filtration.

Therefore, to
show  the equality,  when
$G$ is
$SU_q(2)$, of the algebra in Popa's construction with the fixed point algebra,
  it is therefore sufficient to check that the intersections of these
two algebras with the space
$$B(V_\pi)Q_0 B(V_\pi)... Q_0 B(V_\pi) $$ coincide. (This comes also 
to to the determine the
fixed point algebra in spaces of the form $B(V_\pi)^{\otimes^n}$, by 
formally replacing
elements in $Q_0$ with tensor product sign. Counting dimensions of 
fixed point algebra
could probably give an alternative for the argument bellow).

  Since we allready have shown the reverse equality, it is sufficient 
to check that any
  element in the fixed point algebra intersected with one of the above 
finite dimensional
spaces $B(V_\pi)Q_0 B(V_\pi)... Q_0 B(V_\pi) $ is one of the elemnts 
that respects the open and
closing paranthesis structure of the  reccursive construction of the 
Popa's algebra.
\def\a{\alpha}
\def\b{\beta}

So assume that for given scalars $\mu_{\a_0...\b_n}$  we have an 
fixed point element of
the form
$$(\sum_{\alpha_0,...,\beta_n}\mu_{\alpha_0....\beta_n}
e_{\alpha_0\beta_0}q_1e_{\alpha_1\beta_1}q_2.... 
q_ne_{\alpha_n\beta_n})\otimes 1=
$$
\begin{equation}
\label{desimplificat}
\sum_{r_0,...,s_n}e_{r_0s_0}q_1e_{r_1s_1}...q_ne_{r_ns_n}\otimes
( \sum_{\alpha_0,...,\beta_n}\mu_{\alpha_0....\beta_n} u_{r_0\a_0}u^*_{s_0\b_0}
u_{r_1\a_1}u^*_{s_1\b_1}... u_{r_n\a_n}u^*_{s_n\b_n}).
\end{equation}
Then the fixed point condition comes to the following condition that 
has to hold
for all $r_0,...,s_n$:

\begin{equation}
\label{simplificarea}
\sum_{\alpha_0,...,\beta_n}\mu_{\alpha_0....\beta_n} u_{r_0\a_0}u^*_{s_0\b_0}
u_{r_1\a_1}u^*_{s_1\b_1}... u_{r_n\a_n}u^*_{s_n\b_n}=\mu_{r_0,....,s_n}.
\end{equation}

By denoting $\Theta$ the matrix 
$(\mu_{\alpha_0....\beta_n})_{\alpha_0....\beta_n}$, and by
$U$ the matrix $(u_{ij})_{ij}$ and by $U^\#$ the matrix 
$(u_{ji})_{ij}$, the relation
\ref{simplificarea} comes to the following equation:
\begin{equation}
\label{simplif}
  (U\otimes U^\#....\otimes U\otimes U^\#)\Theta=\Theta.
\end{equation}

Let $I_{(ik)(rs)}$ and $E_{(ik)(rs)}$ be the matrices defined by
$I_{(ik)(rs)}=\delta_{rs}$ and $E_{(ik)(rs)}=\lambda^{-1}_s\delta_rs$. Then
the equations \ref{eqDe.4}, \ref{eqDe.5} correspond to
\begin{equation}
\label{trans1}
  (U\otimes U^\#)I= I; \ \ \ \ I(U^\# \otimes U) =I,
\end{equation}
\begin{equation}
\label{trans2}
  (U\otimes U^\#)E= E; \ \ \ \ E(U^\# \otimes U) =E,
\end{equation}
where $(U\otimes U^\#)$ and $(U^\# \otimes U)$ are considered as the matrices
indexed as follows:
$(U\otimes U^\#)_{((st)(ij) }=u_{si}u^*_{tj}$ and
$(U^\# \otimes U)_{((st)(ij) }=u_{si}u^*_{tj}$. As in \cite{MuPo}, \cite{Wo3},
the relations \ref{trans1},\ref{trans2} are the only relations 
defining $SU_q(2).$
Thus if a matrix $\Theta$ verifies the relation \ref{simplif}, then 
$\Theta$ is obtained
by reccurence (and linear span) from consecutive applications of the relations
\ref{trans1}, \ref{trans2}.

  To have a reduction we therefore have to use the relation $\sum
_j\lambda_j^{-1} u^*_{sj} u_{tj}=\lambda_s^{-1}\delta_{st}$, or 
$\sum_{j}u_{sj}u^*_{tj}
=\delta_{st}.$ (the other two equtions in \ref{eqDe.4}, \ref{eqDe.5} 
involve a summation index
which is not appropriate for the corresponding sum).
This means that the only possibility to have a simplification in
\ref{simplificarea}, by using the defining relations of $SU_q(2)$, is 
if somewhere
in the sum of \ref{simplificarea}, we have that
$\mu_{\alpha_0....\beta_n}$ splits as a sum of elements the form
\begin{equation}
\label{red1}
\lambda^{-1}_{\b_i}\delta_{{\b_i}{a_{i+1}}}
\times \mu'_{\alpha_0,...,\hat{\b_i},\hat{\a}_{i+1},...,\b_n},
\end{equation}
or in the form
\begin{equation}
\label{red2}
\delta_{{\a_i}{b_{i}}}
\times \mu''_{\alpha_0,...,\hat{\a_i},\hat{\b_{i}},...,\b_n}.
\end{equation}
The symbol $\hat{\ }$ corresponding to omission of the corresponding 
symbol. Moreover
after applying such a simplification as in \ref{red1} and 
respectively \ref{red2} in
\ref{desimplificat} we obtain  a sum  in $r_0,....s_n$ involving a 
corresponding
$\lambda_{s_i}\delta_{s_ir_{i+1}}$ or respectively  $\delta_{r_is_i}$.

Repeated application of the above two procedures outlined in the 
equations \ref{red1},
\ref{red2}
  to the sum in \ref{desimplificat}, will give, once we have finished 
the simplification
procedure to the sum reducing \ref{simplificarea} to scalars, it 
follows that the
fixed point element has exactly the open and closing paranthesis structure in
the indices $\a_0,...,\b_n$  that corresponds
to an element in the algebra derived from Popa's construction.

\end{proof}

We can also use the  the model described
  in the paper \cite{Ra3}, where the structure of
the free product algebra $Q\ast M_2(\Bbb C)$
  is described in terms of the ``infinite'' free semicircular element
in \cite {Ra2} (here  $M_2(\Bbb C)$ is endowed
with a state). We obtain an explicit random
matrix (\cite{Vo})  model for the  structure of the subfactor 
$\mathcal N^\lambda(Q)\subseteq
\mathcal M^\lambda (Q)$, when
$Q$ is a free group factor or
$Q=L^\infty([0,1]$.

We describe this model bellow:

\begin{figure}[tbp]
\setlength{\unitlength}{1bp}
\begin{picture}(360,360)(-12,-42)
\put(0,0){\includegraphics[bb=0 0 432 386,height=308.8bp,width=345.6bp]
{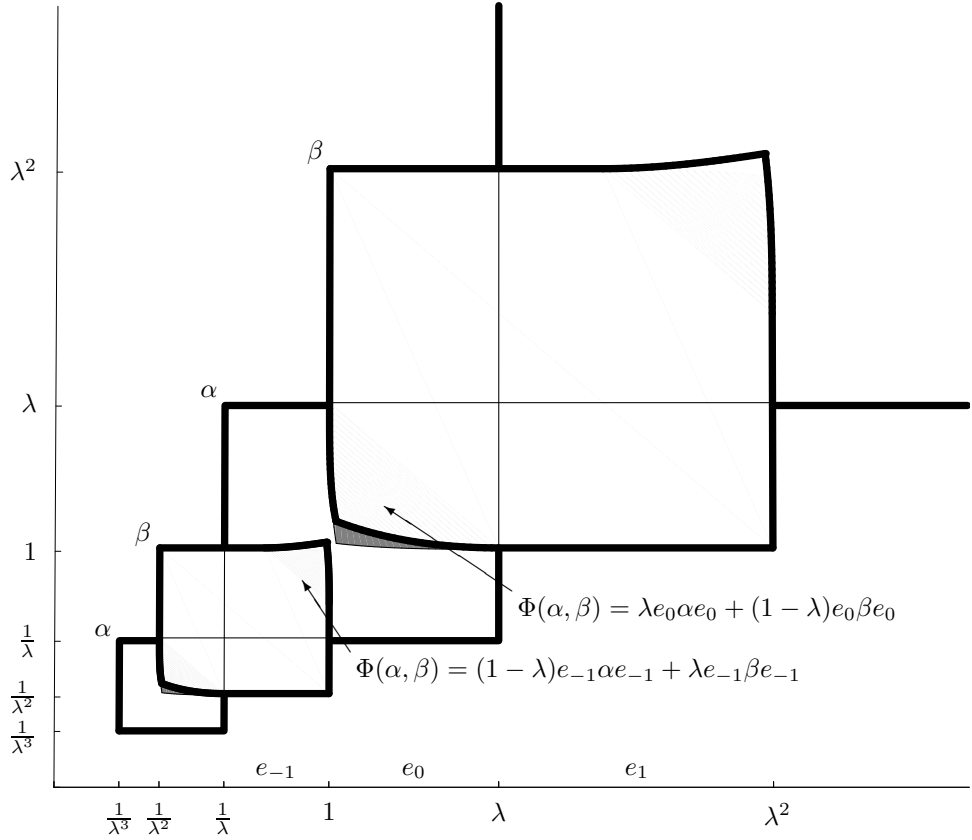}}
\put(25,-6){\makebox(0,0)[t]{$\frac{1}{\lambda^{3}}$}}
\put(40,-6){\makebox(0,0)[t]{$\frac{1}{\lambda^{2}}$}}
\put(64,-6){\makebox(0,0)[t]{$\frac{1}{\lambda}$}}
\put(104,-6){\makebox(0,0)[t]{$1$}}
\put(168,-6){\makebox(0,0)[t]{$\lambda$}}
\put(273,-6){\makebox(0,0)[t]{$\lambda^{2}$}}
\put(-6,19){\makebox(0,0)[r]{$\frac{1}{\lambda^{3}}$}}
\put(-6,34){\makebox(0,0)[r]{$\frac{1}{\lambda^{2}}$}}
\put(-6,55){\makebox(0,0)[r]{$\frac{1}{\lambda}$}}
\put(-6,89){\makebox(0,0)[r]{$1$}}
\put(-6,145){\makebox(0,0)[r]{$\lambda$}}
\put(-6,233){\makebox(0,0)[r]{$\lambda^{2}$}}
\put(22,57){\makebox(0,0)[br]{$\alpha$}}
\put(37,91){\makebox(0,0)[br]{$\beta$}}
\put(62,147){\makebox(0,0)[br]{$\alpha$}}
\put(102,235){\makebox(0,0)[br]{$\beta$}}
\put(112,51){\vector(-2,3){18}}
\put(114,49){\makebox(0,0)[tl]{$\Phi(\alpha,\beta)=(1-
\lambda)e_{-1}\alpha e_{-1}+\lambda e_{-1}\beta e_{-1}$}}
\put(173,74){\vector(-3,2){48}}
\put(175,72){\makebox(0,0)[tl]{$\Phi(\alpha,\beta)=\lambda e_{0}\alpha e_{0}
+(1-\lambda)e_{0}\beta e_{0}$}}
\put(84,9){\makebox(0,0)[t]{$e_{-1}$}}
\put(136,9){\makebox(0,0)[t]{$e_{0}$}}
\put(220,9){\makebox(0,0)[t]{$e_{1}$}}
\end{picture}
\caption{Action of $\Phi$}
\label{Fig1}
\end{figure}

\begin{theorem}
\label{TheoremDe.10}
Let $\lambda$ be a number in $(0,1)$ and let  $D$ be the algebra 
generated by the characteristic functions of the intervals
$e_n=[\lambda^n,\lambda^{n+1}]$, $n\in\Bbb Z$. Let $X,Y$ be two 
``infinite'' free semicircular elements as considered in \cite{Ra2} 
(see also
\cite {DyRa}). Let
$f_i$ be the projections
$e_{2i}+e_{2i+1}$ and let $g_i=e_{2i+1}+e_{2i+2}$. Let $F,G$ be the 
algebras respectively generated by this projections. Let $X_0, Y_0$ 
be the
(bounded) free semicircular family obtained by diagonally carving
$X$ and respectively $Y$ by the projections $f_i$ and respectively 
$g_i$. In \cite{Ra3} it was proved that the algebra $\mathcal A$ 
generated
by $X_0$, $Y_0$ and $D$ is isomorphic to $M_2(\Bbb C)\ast L^{\infty 
}([0,1]$, where on $M_2(\Bbb C)$ we consider a state with weights
$\lambda$ and $1-\lambda$.

Note that  the automorphism of homothety by $\lambda$ (which is a 
restriction of
the one parameter group introduced in \cite{Ra2}),  scales trace by
$\lambda$.

Let $\mathcal B_0, \mathcal B_1$ be the algebras $F'\cap \mathcal A$ 
and respectively
$G'\cap \mathcal A$. Let $E_0, E_1$ be the sum of the even 
(respectively) odd indexed
projections in $(e_i)$.
On the algebra $\mathcal B_0\oplus \mathcal B_1  $ we consider the
following linear map
$$\Phi(\alpha,\beta)=\lambda E_0\alpha E_0+ (1-\lambda) E_0\beta E_0 
+(1-\lambda) E_1\alpha
E_1+
\lambda E_1\beta E_1,
(\alpha,\beta)\in \mathcal B_0\oplus \mathcal B_1.$$
Let $\mathcal M$ be the minimal subalgebra of $\mathcal B_0\oplus 
\mathcal B_1  $ containing $(X_0,Y_0)$ and closed under $\Psi$, and 
let
$\mathcal N$ be the image of $\Phi$.

Then $\mathcal N\subseteq \mathcal M$ is an irreducible subfactor of 
index $\frac {1}{\lambda}+
\frac {1}{1-\lambda}$.

\end{theorem}
\begin {proof}
Indeed via the the identification of $M_2(\Bbb C)\star 
L^\infty([0,1])$ this corresponds exactly to the subfactor 
constructed in
\ref{LemDe.4}, by reducing by the projection $f_0$ the first 
component of $\mathcal C$.
\end{proof}

\section{\label{De2} Determination of the cross product algebra  }

In this section we determine the structure of the crossed product algebra
$(Q\ast B(V_\pi))\rtimes G$.  This, by Takesaki Duality, is in turn 
used to determine, the
structure of the fixed point algebra.

The fixed point algebra could be a type $III$ factor or a type $II_1$ factor.
In the case of the factor in the previous section
  the fixed point algebra is a type
$II_1$ factor which makes it easier to determine the structure of 
this algebra. This
  is similar to the results in \cite{Na}, \cite {Ko}, \cite {KW}, 
where the fixed
point algebra of an infinite tensor product of $B(V_\pi)$ under a 
natural action of
$SU_q(N)$ is determined to be the Temperly-Lieb algebra.

To establish that the cross product is stably
isomorphic to the fixed
point algebra we will need to use instead a quotient group 
``$SU_q(2)\slash \Bbb Z_2=
SO_q(3)$''( see \cite {Pod}). Indeed, in the
  statement of theorem \ref{TheDe.5} the action
  of $SU_q(N)$ on $B(V_\pi)$ is not faithful, in the sense  that its tensor
products do not contain all other representations of the group.

  We will use an idea
from Wasserman's paper (\cite{Wa}) that in the case of an infinite 
tensor product action of
a compact group, by a faithful, selfadjoint, unitary irreducible 
action of a compact group,
the Takesaki's duality gives a stable isomorphism between the fixed 
point algebra
and the cross product.
  We are also indebted to Y. Ueda for
  pointing out to us that part of the statement in theorem 5.6
  in \cite {Na} doesn't hold, as also mentioned in the abstract of the paper
by M.  Izumi \cite{Iz}).

\begin {proposition}
\label{proso3} Let $V$ be the Hilbert space of the fundamental 
representation of $SU_q(2)$.
  The adjoint coaction of $SU_q(2)$ on $B(V)$ corresponds to
a coaction $\alpha$ of $SO_q(3)$ on $B(V)$.
  Moreover the fixed point algebras for any of the two (free product) 
coactions on $Q\ast B(V)$
coincide. Thus
$$(Q\ast B(V))^{SU_q(2)}=(Q\ast B(V))^{SO_q(3)}$$
\end{proposition}
\begin{proof}
Let as in \cite {Pod} $d_0$, $d_{1\slash2}$, $d_1$, ... be the
  representations of $SU_q(2)$, with $d_{1\slash2}$ the fundamental 
representation.
Then the fundamental unitary for $SO_q(3)$ is given by the
  unitary representing $d_1$. As $V\oplus\overline V$ is $d_0\oplus d_{1}$  it
follows that the  coaction of
  $SU_q(2)$ gives the coaction $\alpha$ of $SO_q(3)$ on $B(V)$, (for 
classical groups that
means that the representation comes from a representation of the quotient).
  Since the unitary implementing both coactions is the same, it also 
follows that we
have the same fixed point algebra, (in terms of groups representations
  if we have a representation that factors through the quotient, then we have
the same fixed point algebra).

\end{proof}
In the next proposition we prove that for any finite dimensional
  Hilbert space $V$, if the  tensor product of a selfadjoint  coaction $\alpha$
of a Woronowicz quantum group $G$ on $B(V)$ contains any
other finite dimensional unitary representation of the quantum group,
    then  the coaction  on
$Q\ast B(V)$ is semidual in the sense of \cite{NaTa} \cite{Wa}.
  Thus the fixed point algebra is Morita equivalent to the cross product
by the action. We will use Takesaki Duality, in the sense described 
in \cite{Na}.
\begin{proposition}
\label{ProDe.6}Let $A$ be the function algebra of   Woronowicz a 
compact quantum group $G$, with faithful Haar state. Let $L^2(A)$ be 
the Hilbert
space associated to the Haar measure on
$A$. Assume that
$\alpha$ is a faithful selfadjoint corepresentation of $G$ (i.e a 
corepresentation such that its tensor product contains any
  other unitary, finite dimensional irreducible representation of $G$). Let
$Q$ be a
$II_1$ factor (or a diffuse abelian algebra).   Then
$$(Q\star B(V))\rtimes A\cong (Q\star B(V))^G \otimes B(L^2(A)).$$
\end{proposition}
\begin{proof} Let $\mathcal A= Q\star B(V)$. By the argument in the 
proof of Theorem 5.6 in \cite{Na} (and using also lemma 20 in 
\cite{Bo},
  as pointed out in\cite{ShUe})
it is sufficient to prove that for any irreducible unitary coaction 
$\hat{\alpha}$ of $G$, the corresponding spectral subspace
(\cite{Bo},\cite{Wa},
\cite{HLS})
$\mathcal A_{\hat{\alpha}}$ is non-trivial. But, since $B(V)$ appears 
in a tensor product situation in
$Q\star B(V)$ this follows from the fact that the tensor product of 
the representation $\alpha$ with itself contains any other
representation (see  \cite{Wa}, \cite{NaTa}).
\end{proof}

The following was proved in \cite{Ue1}.  It expresses the natural 
fact that the crossproduct
distributes when we have an action on a free product  of two 
algebras, such that the  coaction
is trivial on one of the two factors in the free product.

  Note that one has to specify a state
on  the algebra over which we amalgamate in order to get a von Neumann algebra.

\begin{proposition}
\label{ProDe.7}(\cite{Ue1}) Let $A$ be
  the function algebra of  a Woronowicz
  quantum group $G$, let $\alpha$ be a corepresentation
of $G$ on
  the bounded linear operators $B(V)$
  acting on a finite dimensional Hilbert space $V$. Let
$\phi_\alpha$ be a $\alpha$ invariant faithful
   state on  $B(V)$, that is such
  that $(\phi_\alpha\otimes Id )(\alpha(x))=\phi_\alpha(x) 1$,
  for all $x$ in $B(V)$.
Let
$\hat A$  be .the dual algebra.
  Then  we have the following isomorphism
  $$(Q\ast B(V))\rtimes A\cong (Q\otimes \hat A)\ast_{\hat A} (B(V)\rtimes A).$$
The amalgamated free product is
  with respect to the canonical conditional expectation of
$B(V)\rtimes A$ onto $\hat A$ (i.e
  the restriction of $\phi_\alpha\otimes Id) $. Moreover  the free 
amalgamated von Neumann
algebra on the right hand side of the above equality is
  determined by endowing the algebra
$\hat A$  with a state (or weight)
  that is the restriction of a faithful normal state on $B(L^2(A))$.
\end {proposition}

We note in the case of $SO_q(3)$ the
  hypothesis are satisfied. Indeed we have the following:

\begin{lemma}
\label{lempi}We use the notations from
  the above lemma. For the action $\alpha$ in \ref{proso3}, if we 
consider the state
$\phi_\alpha=\tau=tr (Q_\pi\cdot)$, where
$Q_\pi$ is the operator associated in \cite
  {Ba} to the  fundamental representation $\pi$ of $SU_q(2)$, then 
$\phi_\alpha$ is
$SO_q(3)$ invariant.
\end{lemma}\begin{proof}. Indeed the unitary
  implementing the adjoint corepresentation $Ad\ \pi$ is the same as 
the unitary implementing
$\alpha$. Since $\tau$ has this property, the same will hold for $\phi_\alpha$.

\end{proof}

The algebras $\hat A\subseteq B(V)\rtimes A$
are discrete   and  the structure of the
  corresponding inclusion matrix can be easily
described in the case when the representation
ring of $A$ is the same as the representation ring of a classical 
compact group.
We are indebted to V. Toledano, D. Bisch,
  M. Pimsner and G.Nagy for pointing us  the
  more precise description of this inclusion matrix,
which essentially appears in \cite{Wa}. We will
  not use the explicit description of this matrix.

\begin{lemma}
\label{ProDe8}Let $A$ be the Woronowicz
  function algebra of a quantum group $G$. Let $\hat A$ be the dual 
algebra and let $\pi$ be a
finite dimensional unitary corepresentation of $A$ on a Hilbert space $V_\pi$.
  Then  $ B(V_\pi)\rtimes A$ is isomorphic
  to $B(V_\pi)\otimes \hat A$. Moreover $\hat A\subseteq 
B(V_\pi)\rtimes A$ has a Bratelli
inclusion matrix  with finite multiplicities.
\end{lemma}
\begin{proof} By definition $ B(V_\pi)\rtimes A$
  is isomorphic to the von Neumann algebra generated in 
$B(V_\pi)\otimes B(L^2(A)$ by
$u^*(B(V(\pi)\otimes 1)u$ and $\hat A$.
But $u(1\otimes \hat A)u^*$ is contained in $B(V_\pi)\otimes \hat A$, since
$Ad\ u$ acting on $B(L^2(A))$ maps the matrix
   coefficients (viewed as elements of $\hat A$) corresponding to a 
finite dimensional unitary
representation
$u^\alpha$ of $A$ into
  a linear combination of
  the matrix coefficients (viewed as elements in $\hat A$) of the 
representation $u\otimes
u^\alpha$.

  We describe this map
in classical setting and then explain the modifications needed for a 
quantum compact group.
Indeed in a classical setting, if $G$ is a compact group, let
$V$ be the fundamental unitary  viewed as an element inin 
$L^\infty(G)\otimes B(L^2(G))$
defined by
$V(g)=\lambda_g$, $g\in G$, where $\lambda_g, g\in G$ is the left 
regular representation
of $G$ into the unitary group on $L^2(G)$.

Assume $\pi(g)=(u_g)\in \mathcal U(B(V_\pi)), g\in G$ is
  a finite dimensional unitary representation of $G$ on a Hilbert 
space $V_\pi$. Let
$U$ be the unitary in $L^\infty(G) \otimes B(V_\pi)$ given by $u$.
Then in $B(L^2(G)) \otimes B(V_\pi)$ the following holds (the 
absorption principle for the left
regular representation):
$$U^* (\lambda_g \otimes 1) U=\lambda_g\otimes u_g.$$
By using the fundamental unitary $V=V_{12}$ defined above, in
$L^\infty(G)\otimes B(L^2(G)) \otimes B(V_\pi)$, the last equality gives
(by using Woronowicz's notations)
\begin{equation}
\label{clasabs}
U^*_{23}V_{12} U_{23}= V_{12}U_{13}.
\end{equation}

Let $e_{ij}$ be a matrix unit for $B(V_\pi)$ and let $u_{ij}(g)$ be 
the matrix coefficients
for $u_g$ in this matrix unit.
For convolution operators in $L(G)$ the equality \ref{clasabs} 
corresponds, for $f$ in $L^1(G)$ to the
following equality in $B(L^2(G))\otimes B(V_\pi)$.
\begin{equation}
\label{clasfus}
U^*(\int f(g)\lambda_g dg) U=
  \int  \sum_{ij}(f(g) u_{ij}(g))(1\otimes e_{ij}) \lambda_g dg.
\end{equation}
This gives, using the decomposition of tensor products of irreducible
  unitary representations of $G$,  a complete description of the
inclusion $L(G)\subseteq B(V_\pi)\rtimes G$ (since we have described
the inclusion $U(L(G)\otimes 1) U^*$ into $B(V_\pi)\otimes L(G)$ in 
terms of fusion rules
of tensor by the representation $\pi$ (on matrix coefficients).

If we want to generalize the above statement for arbitrary quantum 
groups, we will use
the corresponding absorption principle for arbitrary quantum groups 
described in
formula 5.11 in \cite {Wo2}  which replaces \ref{clasabs}.
This formula now holds in $A\otimes B(L^2(A)\otimes B(V_\pi)$ (where 
$\pi$ is the
representation in our statement), $V$ is now one of the fundamental unitary in
  $A\otimes B(L^2(A)$ replacing   the $V$ used for the classical case.
In this context a convolution operator $\int f(g) \lambda_gdg$ is 
replaced (in $A\otimes B(L^2(A))$ by
$$(h(f\cdot) \otimes Id)(V)\in \hat A, f\in A.$$
  This is a generic element in a weakly dense subspace
of $\hat A$. (here we use
$V_{12}$ which is the flip of what is usually the fundamental unitary).
The formula \ref{clasfus} now reads in $A\otimes B(L^2(A))$ as follows
$$U_{23} \lbrack (h(f\cdot) \otimes Id)(V_{12})\rbrack 
U^*_{23}\otimes Id_{B(V_\pi)}=$$
$$(h(f\cdot)\otimes Id\otimes Id_{B(V_\pi)}(U_{23}V_{12}U_{23}^*).$$
By using the formula (5.11) in the paper \cite{Wo2} we get that this 
is further equal
to (and using a matrix unit $e_{ij}$ in $B(V_\pi)$ with respect to which $u \in
A\otimes B(V_\pi)$ has components $u_{ij})$)
$$
(h(f\cdot)\otimes Id\otimes Id_{B(V_\pi)})(V_{12}U_{13})$$ which is 
thus the following
element in $B(V_\pi) \otimes \hat A$,
\begin{equation}
\label{qufus}
\sum_{ij} (h(fu_{ij}\cdot)\otimes Id)(V)\otimes e_{ij}.
\end{equation}

This shows that in $A\otimes B(L^2(A)$,
  the transformation  $Ad\ u$ maps $ B(V_\pi)\rtimes A$
onto $B(V_\pi)\otimes \hat A$.
\end{proof}
\begin {remark}
\label{refus}
If $G=SU_q(2)$, and $\pi$ is the fundamental representation, since 
the representation ring of
$G$ is the same as the classical one, it follows that the inclusion 
matrix of $\hat A$ into
$B(V_\pi)\rtimes A$, is the same a sinthe classical one, which
is  a matrix of type $A_\infty$.
\end{remark}
We need to apply the previous lemma to the
  case of the coaction $\alpha$ of $SO_q(3)$ on $B(V)$ which was described in
\ref{proso3}.
\begin{lemma}
\label{lemmacuB$}
  Let $\hat B$ be the dual Woronowicz
  algebra for $SO_q(3)$.
  With the notation in proposition \ref{proso3}
  we have that $B(V)\rtimes SO_q(3)$ is isomorphic to $B(V)\otimes 
\hat B$. Moreover the
inclusion matrix $\hat B\subseteq B(V)\rtimes SO_q(3)$ has a Bratelli 
inclusion matrix with
finite multiplicities.
\end{lemma}

\begin{proof}Denote $B=SO_q(3)$ and $A=SU_q(2)$ and
  let $\hat B,\hat A$ be the dual algebras.
The representation $\alpha: B(V)\rightarrow B(V) \otimes B$
  is the restriction of the adjoint corepresentation induced
  by $\pi$. This holds because
$\alpha(x)=u^*(x\otimes 1 )u$ belongs to $B(V)\otimes B$,
since $B$ is generated by the matrix coefficients of $d_1$.
But then if $P$ is the projection (in $B(L^2(A)$) onto $L^2(B)$, then
$$B(V)\rtimes B=P(B(V)\rtimes A)P,$$
since $P\hat AP=\hat B$ and since $B$ acts on $L^2(A)$ as the Haar 
measure on $B$ is the
restriction of Haar measure on $A$ (\cite{Pod}).
\end{proof}

The following lemma is a direct consequence of the method used in 
\cite{Ra}. There it was proved that amalgamated free products of the 
type
$(L(F_N)\otimes D)\ast_D C$, where $D\subseteq C$ is an inclusion of 
discrete von Neumann algebras, $C$ with a faithful trace, are
  isomorphic to
a free group factors.
\begin{lemma}
\label{LeDe.9} Let $D\subseteq C$ be   von Neumann algebras that are 
infinite sums of algebras of finite dimensional matrices, and such 
that the
Bratelli inclusion matrix has finite multiplicities and
$\mathcal Z(C)\cap\mathcal Z(D)=\Bbb C1$.  Let $Q$ be a free group 
factor or a diffuse abelian
von Neumann algebra. Assume that $D$ is endowed with a   (semi)finite 
faithful trace $tr$  and that the amalgamation is
performed (\cite{Vo}, \cite {Po1}) with respect to a
normal faithful  conditional expectation from $C$ onto $D$, which 
then gives a trace state on the amalgamated free product,
(which by GNS construction gives the amalgamated free product von 
Neumann algebra).

Assume that the amalgamated free product von Neumann algebra is 
finite or semifinite.

Then $tr\circ E$ is a faithful  trace  on $C$ and
$\mathcal L(F_N) \otimes D
\star_D C$ is Morita equivalent to  a free group factor.
\end{lemma}
\begin {proof}
Indeed if $\tau \circ E$ is not a trace then the modular group of 
$\tau\circ E$ will be non-trivial on off diagonal elements. If the 
trace
is finite this is exactly the content of  theorem 5.1 in \cite{Ra} 
(see also \cite {Shly} for a different, more recent proof). If the 
trace is
infinite, the arguments in the proof of the theorem   mentioned 
above can  obviously be modified to handle this case (by changing the 
principle of
counting the blocks as in the construction of a one parameter group 
of automorphisms in \cite{Ra2}).
\end{proof}

\begin{remark}
\label{type2} With the notations in theorem \ref{TheDe.5}, one can 
prove directly that factor
$(Q\ast B(V_\pi))\rtimes SU_q(2)$ is  a finite  type II von Neumann algebra.
Indeed because of  lemma \ref{ProDe.7}, it is sufficient to that the inclusion
$\hat A\subseteq B(V_\pi)\rtimes A$, and the corresponding 
conditional expectation
$tr(\phi_\alpha\cdot )\otimes Id$ (see also lemma  \ref{lempi})
  from
$B(V_\pi)\rtimes A$ onto $A$, has the property that the composition 
of the trace
on $\hat A$ with the conditional expectation is again a trace.
The trace on $\hat A$ that we are considering here is the restriction 
to $\hat A \subseteq
B(L^2(A)$ of the canonical trace on $B(L^2(A))$.

  Let $p_\pi$ be the projection from $\hat A$
onto
$B(V_\pi)$ (by viewing $\hat A$ as the direct sum over all bounded 
linear operators on
  the Hilbert spaces  of an enumaration of the irreducible finite 
dimensional unitary
representations of $A$,
(\cite {KuNa},
\cite{PoWo})).

Let $\hat \phi :\hat A\rightarrow \hat A\otimes \hat A$ be dual 
comultiplication map.
Then the inclusion $\hat A\subseteq B(V_\pi)\rtimes A$, (in the identification
$B(V_\pi)\rtimes A\cong B(V_\pi)\otimes \hat A$, described in 
Proposition \ref {ProDe8})
  is, (because of the description in \cite{PoWo}) exactly the map
$$ (p_\pi\otimes Id)\hat\phi: \hat A\rightarrow B(V_\pi)\otimes \hat A.$$
  Moreover if
$\hat \gamma$ is the density matrix in $\hat A$ of the dual Haar 
measure $\hat h$  on $\hat A$,
then $\hat \phi(\hat \gamma)=\hat \gamma\otimes \hat \gamma$ (\cite 
{Ya}). But this
amount exactly to the fact that the composition of the conditional expectation
$tr(\phi_\alpha\cdot )\otimes Id$ with $\hat h$ gives the state 
$tr(\phi_\alpha\cdot )\otimes
\hat h$ on $B(V_\pi)\rtimes A\cong B(V_\pi)\otimes \hat A$. Hence, 
the composition
of the trace on $\hat A$ (coming from $B(L^2(A)$) with the 
conditional expectation is a trace.

\end{remark}

\begin{corollary}
\label{CorDe.8} The factor $\mathcal M^\lambda(\mathcal L(F_k))$ in 
Theorem \ref{TheDe.5} (for $N=2$ ) is isomorphic to $\mathcal 
L(F_\infty)$.
\end {corollary}

\begin{proof} From \ref{LeDe.9}, \ref{ProDe.7} it follows that 
$(Q\ast B(V_\pi)\rtimes SO_q(3)$ is isomorphic to
$\mathcal L(F_\infty)\otimes B(H)$, where $H$ is an infinite 
dimensional space. The result  follows  now from
\ref{TheDe.5}, \ref{proso3}, since $d_1$ is a faithful representation 
of $SO_q(3)$ (no adjoint is required here).

\end{proof}

In this way we reobtain, by a different method, the result that was 
recently proved by D. Shlyakhtenko and Y. Ueda in \cite{ShUe}.
We do not know
if the subfactors obtained by using the method in this paper, which 
are derived from Popa's construction of irreducible subfactors from 
the
Temperly-Lieb algebra   coincide with those constructed in the  paper
\cite{ShUe}. Both of the two subfactors  have $A_\infty$ invariants and both
are  the fixed point algebra of a coaction of a quantum group, 
although in one case the fixed point algebra is a type III factor 
while in the
present
  case this is a type
$II_1$ factor.
\begin{corollary}
In
particular
$\mathcal L(F_\infty)$ has irreducible subfactors of index $\lambda$ 
(and type $A_\infty$) for all index values bigger than 4.
\label{CorDe.9}

\end{corollary}
The previous result, since we are using a non-tracial version of \cite{Po1} is
strong evidence to Popa's
  conjecture that the subfactors
constructed in  the  breakthrough papers \cite{Po1} (or 
\cite{Po2},\cite{Po3}) are isomorphic
to free group factors . The only case in which
Popa's subfactors (\cite{Po1},\cite{Po2}) are  known to be free group 
factors is for index
values less than 4 (\cite {Ra}).

\bibliographystyle{bftalpha}

\ifx\undefined\bysame

\newcommand{\bysame}{\leavevmode\hbox to3em{\hrulefill}\,}

\fi


\begin{thebibliography}{CoBull}


\bibitem[BaSa]{BaSa}

  Baaj, Saad; Skandalis, Georges
{\em Unitaires multiplicatifs et dualite pour les produits croises de 
$C\sp*$-algebres.
}
Ann. Sci. Ec. Norm. Super., IV. Ser. {\bf 26}, No.4, 452-488 (1993).


\bibitem[Ba]{Ba}
  Teodor Banica
Representations of compact quantum groups and subfactors. (English)
J. Reine Angew. Math. {\bf 509}, 167-198 (1999).

\bibitem[Ba1]{Ba3}
  Banica, Teodor
{\em Theorie des representations du groupe quantique compact libre $O(n)$. }

C. R. Acad. Sci., Paris, Ser. I {\bf 322}, No.3, 241-244 (1996).

\bibitem[Bar]{Bar}
  Barnett, Lance
Free product von Neumann algebras of type III. (English)
Proc. Am. Math. Soc. {\bf 123}, No.2, 543-553 (1995).


\bibitem[BlDy]{BlDy}
Etienne Blanchard, Ken Dykema,
{\em  Embeddings of reduced free products of operator algebras}, OA/9911012,
  to appear in Pacific J. Math. 15 pages.

\bibitem[Bo]{Bo}
Boca, Florin P.
{\em Ergodic actions of compact matrix pseudogroups on $C\sp*$-algebras.} in
  Connes, A. (ed.), Recent advances in operator algebras. Collection 
of talks given in the
conference on operator algebras held in Orleans, France in July 1992. 
Paris: Societe
Mathematique de France, Asterisque. {\bf 232}, 93-109 (1995).

\bibitem[Bo1]{Bo1}
Boca, Florin P., Thesis, UCLA, 1992.

\bibitem[Dy]{Dy}
Kenneth~J. Dykema, {\em Amalgamated free products of multi-matrix algebras and
   a construction of subfactors of a free group factor}, Amer. J. Math. {\bf
   117} (1995), no.~6, 1555--1602.

\bibitem[Dy1]{Dy1} Dykema, Kenneth J.
{\em Amalgamated free products of multi-matrix algebras and a 
construction of subfactors of a
free group factor} Am. J. Math. {\bf 117}, No.6, 1555-1602 (1995).

\bibitem[DyRa]{DyRa}
Dykema, Kenneth J.; Radulescu, Florin
{\em Compressions of free products of von Neumann algebras}
Math. Ann. {\bf 316}, No.1, 61-82 (2000).



\bibitem[Iz]{Iz}
  Izumi, Masaki
{\em Actions of compact quantum groups on operator algebras.}
RIMS Kokyuroku {\bf 1024}, 55-60 (1998)

\bibitem[Jo]{Jo}
V.~F.~R. Jones, {\em Index for subfactors}, Invent. Math. {\bf 72} (1983),
   no.~1, 1--25.


\bibitem[Jo1]{Jo1}
Jones, Vaughan F.R.
{\em Subfactors and knots}. Expository lectures from the CBMS 
regional conference, held at the
US Naval Academy, Annapolis, USA, June 5-11, 1988.
Regional Conference Series in Mathematics. {\bf 80}. Providence, RI: 
American Mathematical
Society. viii, 113 p. (1991).




\bibitem[HLS]{HLS}
  Hoegh-Krohn, R.; Landstad, M.B.; Stormer, E.
{\em Compact ergodic groups of automorphisms},
Ann. Math., II. Ser. {\bf 114}, 75-86 (1981).

\bibitem[Ko]{Ko}
  Konishi, Yuji,
{\em A note on actions of compact matrix groups on von Neumann algebra},
Hihonkai Math. J. {\bf 3},, 23-29, (1992).

\bibitem[KNW]{KW}
  Konishi, Yuji; Nagisa, Masaru; Watatani, Yasuo,
{\em Some remarks on actions of compact matrix quantum groups on 
$C\sp \ast$- algebras.}
  Pac. J. Math. {\bf 153}, No.1, 119-128 (1992).
\bibitem [KuNa]{KuNa}
Kurose, Hideki; Nakagami, Yoshiomi,
{\em
Compact Hopf *-algebras, quantum enveloping algebras and dual 
Woronowicz algebras
  for quantum Lorentz groups}
Int. J. Math. {\bf 8}, No.7, 959-997 (1997).

\bibitem[MaNa]{MaNa}
  Masuda, Tetsuya; Nakagami, Yoshiomi
{\em A von Neumann algebra framework for the duality of the quantum groups,}
Publ. Res. Inst. Math. Sci. {\bf 30}, No.5, 799-850 (1994).


\bibitem[Nag]{Nag}
  Nagy, Gabriel
On the Haar measure of the quantum $SU(N)$ group. (English)
Commun. Math. Phys. 153, No.2, 217-228 (1993). MSC 1991: *17B37 
46G12, Reviewer: N.Andruskiewitsch


\bibitem[NaTa]{NaTa}
  Nakagami, Yoshiomi; Takesaki, Masamichi
{\em Duality for crossed products of von Neumann algebras.}
Lecture Notes in Mathematics. {\bf 731}. Berlin - Heidelberg - New York:
Springer-Verlag. IX,
139 p.  (1979).

\bibitem[Na]{Na}
  Nakagami, Yoshiomi
{\em Takesaki duality for the crossed product by quantum groups.}, in
Araki, Huzihiro (ed.) et al., Quantum and non-commutative analysis.
  Past, present and future perspectives. Proceedings of the 
international Oji seminar on
quantum analysis, Kyoto, Japan, June 25-29, 1992 and the symposium on 
non-commutative
analysis, Kyoto, Japan, June 29-July 2, 1992 dedicated to Prof. 
Huzihiro Araki on the occasion
of his 60th birthday. Dordrecht: Kluwer Academic Publishers. Math. 
Phys. Stud. {\bf 16},
263-281 (1993).

\bibitem[Pod]{Pod}
  Podles, Piotr
{\em Symmetries of quantum spaces. Subgroups and quotient spaces of 
quantum $SU(2)$ and $SO(3)$
groups} Commun. Math. Phys. {\bf 170}, No.1, 1-20 (1995).


\bibitem[PoMu]{MuPo}
  Podles, P.; Mueller, E.
{\em Introduction to quantum groups}
Rev. Math. Phys. {\bf 10}  No.4, 511-551 (1998).

\bibitem[PoWo]{PoWo}
Podles, P.; Woronowicz, S.L.
{\em Quantum deformation of Lorentz group}
Commun. Math. Phys. {\bf 130}, No.2, 381-431 (1990).

\bibitem[Po1]{Po1}
Sorin Popa, {\em Markov traces on universal {J}ones algebras and subfactors of
   finite index}, Invent. Math. {\bf 111} (1993), no.~2, 375--405.



\bibitem[Po2]{Po2}
Sorin Popa, {\em An axiomatization of the lattice of higher relative commutants
   of a subfactor}, Invent. Math. {\bf 120} (1995), no.~3, 427--445.



\bibitem[Po3]{Po3}  Sorin Popa, {\em The universal algebra of a subfactor},
  preprint.

\bibitem[Po4]{Po4}
  Popa, Sorin
{\em Classification of amenable subfactors of type II.}
Acta Math. {\bf 172}, No.2, 163-255 (1994).

\bibitem[PoWa]{PoWa}

  Popa, Sorin; Wassermann, Antony
{\em Actions of compact Lie groups on von Neumann algebras.}
C. R. Acad. Sci., Paris, Ser. I {\bf 315}, No.4, 421-426 (1992).

\bibitem[Ra]{Ra}
Florin R{\u{a}}dulescu, {\em Random matrices, amalgamated free products and
   subfactors of the von {N}eumann algebra of a free group, of noninteger
   index}, Invent. Math. {\bf 115} (1994), no.~2, 347--389.


\bibitem[Ra2]{Ra2}
Florin R{\u{a}}dulescu,
{\em A one parameter group of automorphisms
of ${\mathcal L}(F_ \infty)\otimes B(H)$ scaling
the trace.}
C. R. Acad. Sci., Paris, Ser. I {\bf 314}, No.13, 1027-1032 (1992).


\bibitem[Ra1]{Ra1}
Florin R{\u{a}}dulescu, {\em An invariant for subfactors in the von {N}eumann
   algebra of a free group}, C. R. Acad. Sci. Paris S\'er. I Math. {\bf 316}
   (1993), no.~10, 983--988.




\bibitem[Ra3]{Ra3}
Florin R{\u{a}}dulescu,
{\em A type $\text {III}_ \lambda $ factor with
  core isomorphic to the von Neumann algebra of a free group, tensor 
$B(H)$,} in Connes, A.
(ed.), Recent advances in
  operator algebras. Collection of talks given in the conference on
  operator algebras held in
Orleans, France in July 1992. Paris: Societe Mathematique de France,
Asterisque. {\bf 232}, 203-209
(1995).


\bibitem[Ra4]{Ra4}
Florin R{\u{a}}dulescu,
{\em Subfactors inside free group factors}, in preparation.


\bibitem[Shly]{Shly}
Dimitri Shlyakhtenko, {\em Some applications of
  freeness with amalgamation}, J.
   Reine Angew. Math. {\bf 500} (1998), 191--212.

\bibitem[Shly1]{Shly1}
Shlyakhtenko, Dimitri
{\em Free quasi-free states.}
Pac. J. Math. {\bf 177}, No.2, 329-368 (1997).

\bibitem[ShUe]{ShUe}
Dimitri Shlyakhtenko, Yoshimichi Ueda
{\em Irreducible subfactors of $L(\mathbb F_\infty)$ of index $\lambda>4$
}

MSRI Preprint,
Report number: MSRI {\bf 2000-030 } (2000).








\bibitem[Ue1] {Ue1}
Ueda, Yoshimichi,
{\em On the fixed-point
algebra under a minimal free product-type action of the quantum group
$SU_q(2)$},  Int. Math. Res. Not. 2000, No. {\bf 1}, 35-56 (2000)


\bibitem [Ue2]{Ue2}
Ueda, Yoshimichi
{\em A minimal action of the compact quantum group $SU_q(n)$ on a full factor}
J. Math. Soc. Japan {\bf 51}, No.2, 449-461 (1999).

\bibitem [Ue3]{Ue3}
Ueda, Yoshimichi
{\em Amalgamated free product over Cartan subalgebra}
Pac. J. Math. {\bf 191}, No.2, 359-392 (1999).


\bibitem[Vo]{Vo}
Dan Voiculescu, {\em Limit laws for random matrices and free products}, Invent.
   Math. {\bf 104} (1991), no.~1, 201--220.



\bibitem[Wa]{Wa}
  Wassermann, Antony,
{\em Ergodic actions of compact groups on operator algebras.
  I: General theory.}
Ann. Math., II. Ser. {\bf 130}, No.~2, 273-319 (1989).


\bibitem[We]{We}
Hans Wenzl, {\em Hecke algebras of
  type ${A}\sb n$ and subfactors}, Invent.
   Math. {\bf 92} (1988), no.~2, 349--383.







\bibitem[Wo]{Wo}
Woronowicz, S.L.
{\em Compact matrix pseudogroups},
Commun. Math. Phys. {\bf 111}, 613-665 (1987). MSC 1991: *58Z05 58H05


\bibitem[Wo1]{Wo1}
  Woronowicz, S.L.
{\em Tannaka-Krein duality for compact matrix pseudogroups. Twisted 
SU(N) groups}.
Invent. Math.  {\bf 93}, No.1, 35-76 (1988).

\bibitem[Wo2]{Wo2}
  Woronowicz, S.L.
{\em Compact quantum groups,} in
Connes, A. (ed.) et al., Quantum symmetries/ Symetries quantiques.
  Proceedings of the Les Houches summer school, Session LXIV, Les 
Houches, France, August
1 - September 8, 1995. Amsterdam: North-Holland. 845-884 (1998).

\bibitem[Wo3]{Wo3}
  Woronowicz, S.L.
{\em Tannaka-Krein duality for compact matrix pseudogroups. Twisted 
SU(N) groups,}
Invent. Math. {\bf 93}, No.1, 35-76 (1988)

\bibitem[Ya]{Ya}
Yamagami, Shigeru
{\em On unitary representation theories of compact quantum groups}
  Commun. Math. Phys. {\bf 167}, No.3, 509-529 (1995).



\end{thebibliography}
\end{document}